# INNOVATING MATHEMATICS
# BY RESEARCH APPROACHES IN INDONESIA
# FROM MIDDLE SCHOOLS TO UNDERGRADUATES

Hanna Arini Parhusip
Mathematics Departement –Science and Mathematics Faculty
Universitas Kristen Satya Wacana
www.uksw.edu

**Abstract.** Innovating mathematics by research approaches in Indonesia from middle schools to undergraduates are explained here where geometry is typical example to be innovated. Topics in plane geometry to be platonic solids and non platonic solids such as cylinder and spheres are innovated to be new surfaces and media. Topics in calculus and complex functions particularly parametric curves are used to be new curves for curves stitching, motifs in several kinds of souvenirs and accessories. Additionally, local culture is introduced through batik painting activity using algebraic surfaces drawn with Surfer. Leonardo Dome and Leonardo Bridge with Rinus Reolofs like are introduced to attract students for mathematical sports activities.

1. **INTRODUCTION**

When I come to the word innovation in mathematics, it turns out that an innovation is a new method or branch in mathematics itself since ancient time. One has proposed the top ten of innovation in mathematics enriching branches in mathematics such as Logarithms, Matrix algebra, Complex numbers, Non-Euclidean geometry, Binary logic Decimal fractions, Zero and Negative numbers Calculus and Arabic numerals. On the hand, an innovation in my point of view is doing a new creativity from the classical part of a special issue in mathematics and treating it to be new performance or new application with integrating other approaches. Additionally, since technologies are growing in this $21^{st}$ century, rapid innovations in mathematics vary heterogeneously depending on the level of a society to adjust the technologies such that trends on innovating mathematics may differ for each different place or country.

Current trends on innovating mathematics in $21^{st}$ century is merely its application in science (Web1) such that curriculum must be redesign (Mahajan,2014) containing the knowledge, skills and character. In the mathematical knowledge, one must know traditional and modern knowledge requiring multidiscipline. In the skills ability, students are required to know how to use leading to creativity, critical thinking, communication and how to collaborate. Finally, a student must learn how to engage in the world with mindfulness, curiosity, courage, resilience, ethics and leadership. Thus, an innovation on mathematics will vary due to many aspects above. Several approaches will be recalled here to present for the readers how to design creativity in some topics in curriculum and students are engaged with their own innovations and creativities.

Innovating mathematics in this paper addresses on how presenting mathematics differently compared to regular meeting in a classroom leading to some interactions among students, educators and other necessary agents to collaborate with. As a result, new materials may be necessarily prepared. The activity involves not only a single topic in a curriculum but also integration of several topics that may appear due to an innovative sense.

## 2. RESEARCH APPROACH FOR OLD AND MODERN GEOMETRICAL TOPICS

There exists no standard theoretical parts for innovating mathematics since the innovations depend on which basic knowledge will be innovated, i.e. the regular formulation and ideas in mathematical topics are explored into new performance, new beneficial and new audience. To gain a maximal innovation, an integration of several topics is necessary. The contents below are mostly few examples where middle or high schools and undergraduate mathematical topics are studied and innovated by research approaches.

### 2.1 Research approach in geometry of platonic solid

In traditional approach, one must learn some basic knowledge in geometry, e.g. definitions of several distances (between a point and line, point to a plane, a line to a plane, and between 2 planes) where this knowledge is considerably neglected by most students. In fact, students could not make the definitions on their owns though the situations are very regular in daily life. The attitude for delivering mathematical definition without communicating it needs leads to no critical thinking and ability to communicate. Instead, students are posed several new creation and innovation with geometrical old patterns as illustrated in Fig 1. Students learn about platonic solids and the properties leading to an innovation of dress created by students which make students proud to promote.

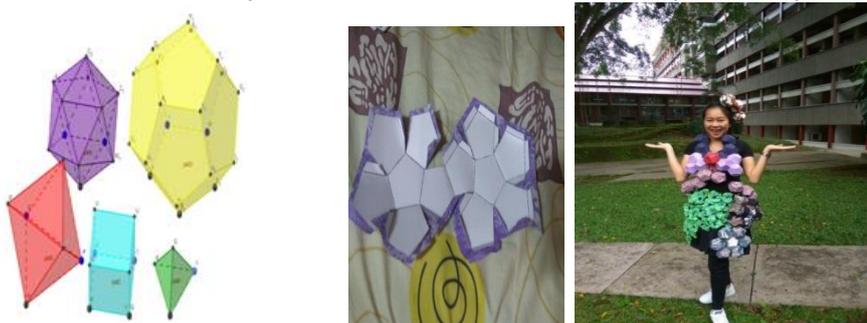

**Figure 1**. Platonic solids in traditional view (left) and in the new performance (middle and right) created by undergraduate mathematics education students in UKSW Indonesia , Sep 2016.

Additional approach can be done by using available software for drawing such as Geogebra since students may show several mathematical properties for each solid.

Rinus Roelofs has more prominent innovations with platonic solids and more innovated geometries as a sculpturist. Roelof's intriguing sculptural objects are made with all kinds of materials including in addition to paper wood metal and acrylic (Web 2). One may find by search engines to overview many kinds of his innovations in mathematical designs to approach a geometry lesson into an interesting topic. Fig.2 depicts several elevated polyhedral of Rinus Roelofs like.

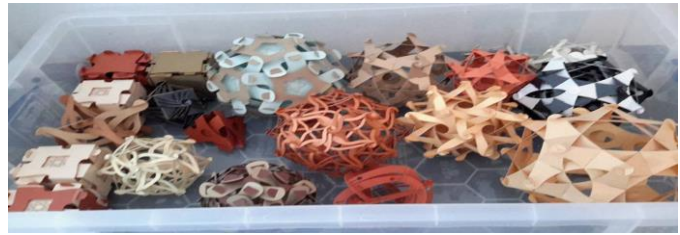

**Figure 2.** Elevated polyhedral created by Rinus Reolofs (source : private photo , taken from IC16, Berlin, Germany, July 2016).

## 2.2 Research approach in Calculus and Complex Functions

Mathematical functions in calculus and complex functions normally are introduced traditionally, i.e. dealing with computation and technical standard for drawing functions. Though introducing functions can also be visualized with software (e.g. Geogebra, MAPLE), the benefits or learning those topics are too formal. Instead, special topics in calculus and complex functions can be combined into a creative activity and the example is shown here.

In high schools level, students may create the equations with excel and innovate the visualization into objects. As stated before, the innovation may create a similarity with nature. Let us consider the parametric equations read as (Parhusip,2016)

$$x = (\sin^2(4\theta) + \cos(4\theta))\cos\theta$$
$$y = (\sin^2(4\theta) + \cos(4\theta))\sin\theta$$

which are obtained from polar equation, i.e.

$$r = \sin^2(4\theta) + \cos(4\theta). \tag{1}$$

The obtained plane curve is similar to a flower called *Icora Javanica* and visualized in Fig.*3* Mapping with complex function  *f(z)*=exp(*z*), puzzle and lampshade are designed and depicted in Fig. 3.

A parametric curve called hypocycloid curve has been modified in various types by considering the resulting curves to be domains for complex functions (Suryaningsih,et.al.,2013) (Parhusip,2014)(Purwoto, et. al.2014) (Parhusip,2015). Results on souvenirs, ornaments, motifs in textiles in collaboration with home industries are produced including some similarity objects in natures, e.g. 2 ornaments are shown in Fig. 4 according to the obtained design as innovative activities in calculus and complex function.

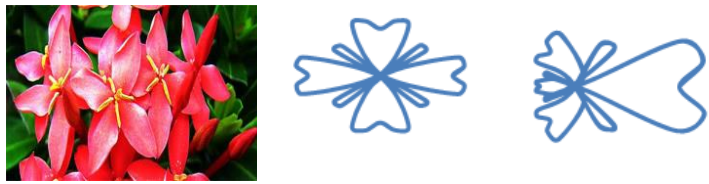

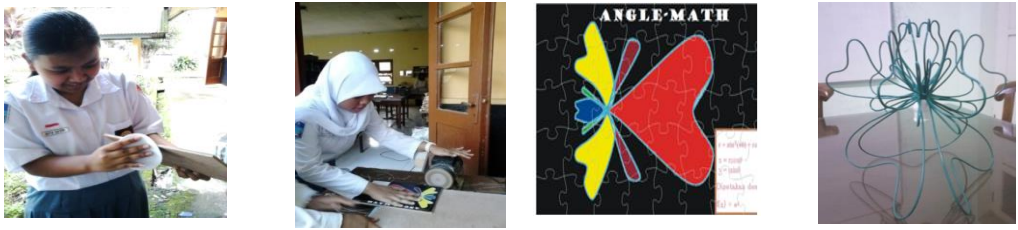

**Figure 3.** High school students activity for creating puzzle and lampshade (Parhusip,2016) from a pair of parametric equation (Eq.1) mapped by *f*(*z*)= exp(*z*).

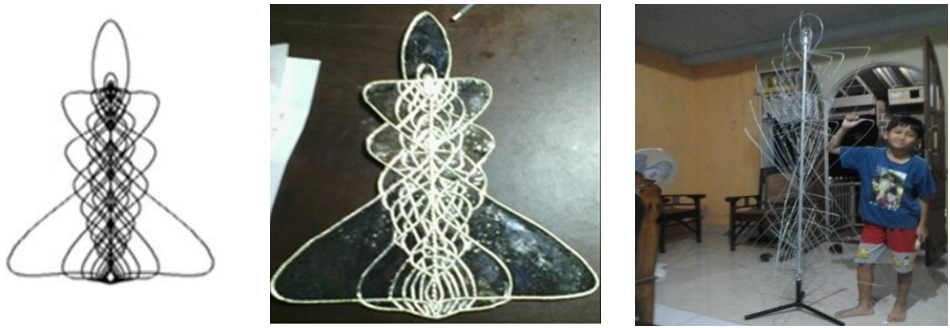

**Figure 4.** Left(above): design from a parametric curve mapped by a complex function (called Rocked Math) ; (left, below) ornament as a souvenir from a silver and a garden decoration (right) (Parhusip,2014)(Parhusip,2015).

### 2.3 The used of functions for data modeling

Typical problem arises from mathematics students is how to deal with data using calculus. Most people refer to employ statistics for data analysis. Unfortunately, some data need to be modeled into mathematical expression before to be analyzed. In classroom, students mostly are introduced functions without any data in hands. By research approach, students are encouraged to model data into an appropriate function leading to optimization problems. Several data from local companies and laboratories in Indonesia have been collected and discussed to give knowledge for beginners in this direction (Parhusip,2016). This approach will not be explored into detail here.

Finally, the method for transferring creatives involve in curriculum must also be taken into account since a standard curriculum requires rigidity in contents. In the section of result and discussion, the activities and the results are explored to inspire readers.

### 3. METHODOLOGY

**3.1 Plane geometry into innovated platonic solids**

We have known that regular plane geometry is normally introduced for students in middle schools until undergraduate, e.g. triangle, rectangle and square, other regular polygon and including circle. Only regular triangle, square and pentagon can be defined into platonic solid yielding to only 5 kinds of platonic solid, i.e. tetrahedron, octahedron, icosahedrons, cube and dodecahedron. The first 3 platonic solids characterized by 3, 4, 5 triangles meet at each vertex respectively where as a cube by 3 squares and dodecahedron by 3 pentagons meet at each vertex. By learning angles at each vertex must be less than $360^0$, one may obtain only 5 platonic solids.

The innovation is devoted to make a special dress with dodecahedron patterns. Drawing and repeating for several patterns of dodecahedron with papers and sticking textiles on these patterns as depicted in Fig.5. Finally, a number of dodecahedrons are used.

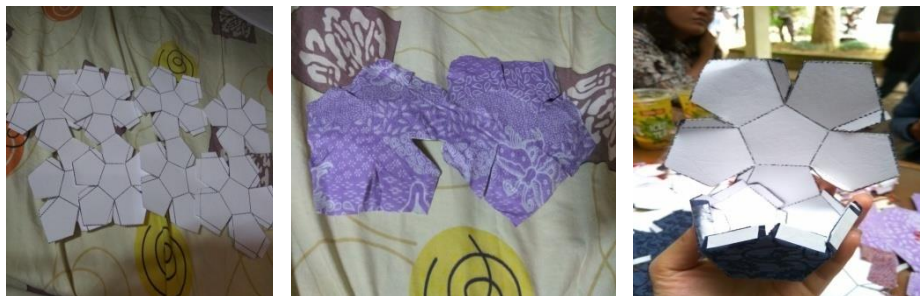

**Figure 5.** Patterns of dodecahedron for special dress

**3.2 Plane geometry for Leonardo Dome and Leonardo Bridge**

Though only rectangles, we may have several activities for developing materials for students to play with geometry. For instance, using bamboos as local materials in Indonesia to create educative activities referring to Rinus Reolofs. Two typical examples here are Leonardo Dome and Leonardo Bridge. The construction of both models are unique,i.e. without using any glue, one may build the models. Additionally physical laws can be explained through these educative materials.

### 4 RESULTS AND DISCUSSION

**4.1 Geometrics activity for middle school students and the related mathematical topics**

As stated above, students must engage into an activity to involve the curriculum. In this level, students should learn about plane geometry including traditional solid objects such as a cylinder, sphere. Students have been introduced with several traditional steps :
the basic knowledge with related topic, i.e. a triangle, a rectangle, a circle, a cylinder, and a sphere and computation for each geometrical object with new approach. In both activities, students are involved in their own discoveries guided by mentors. Some photos of these activities are shown in Fig.4. One of examples of the material is known as Tanggram.

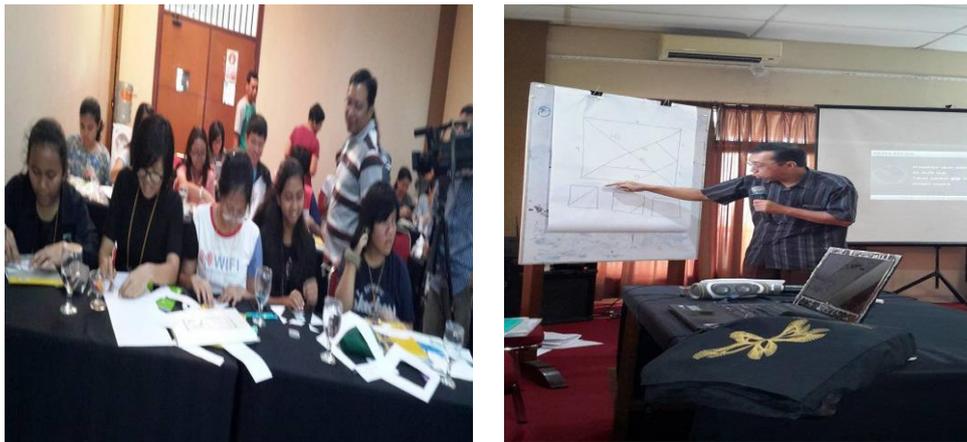

**Figure 4.** Transferring basic knowledge of plane geometry by involving students to create their own materials (Source :private photo, students are from middle schools of Santa Ursula BSD Jakarta and BPK Penabur Cirebon Indonesia under the project Junior Summer Camp, in 10-19 June 2017, Salatiga, Indonesia).

The innovation is then dealt with some local materials leading to other special keywords in mathematics and physics. Software called Surfer is used to introduce algebraic surfaces and the innovation. Cylinder, ellipsoid and ball are typical algebraic surfaces known by middle or high schools students. However combinations of those surfaces have not been discussed and difficult to be illustrated manually. Using a software, some possible other surfaces are easily introduced without learning complete related mathematical terms. Therefore students are introduced to work with Surfer and applied into local culture activity, i.e. create batik motifs according to the drawing results and partly the results are depicted in Fig 4 (Parhusip,2017)

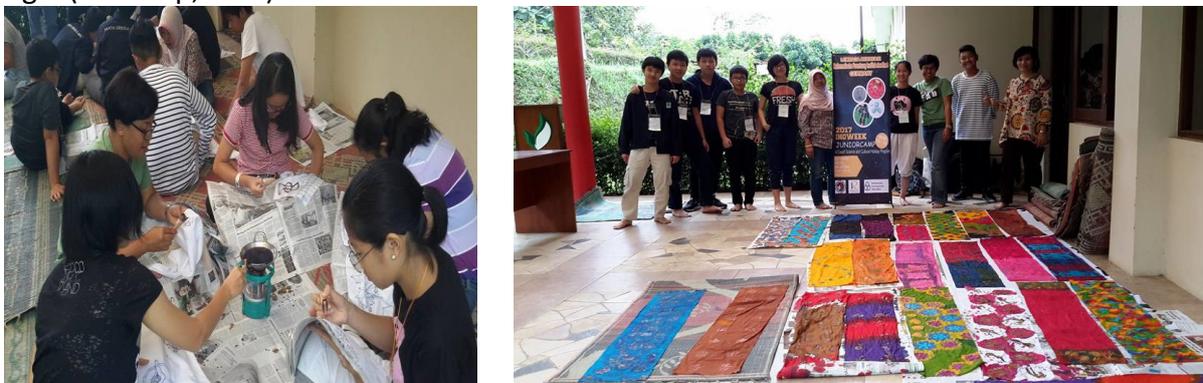

Figure 4. Students are creating and painting motifs from Surfer (left) and the related batik (right) (Source :private photo, students are from middle school of Santa Ursula BSD Jakarta and BPK Penabur Cirebon Indonesia under the project Junior Summer Camp, in 10-19 June 2017,Salatiga, Indonesia June 2017)

Geometrical innovation is then developed with the material adopted from Rinus Reolofs, i.e. building Leonardo Dome and Leonardo Bridge from woods. Before this activity is presented, an architecture need in future could be introduced such that students are motivated to involve.

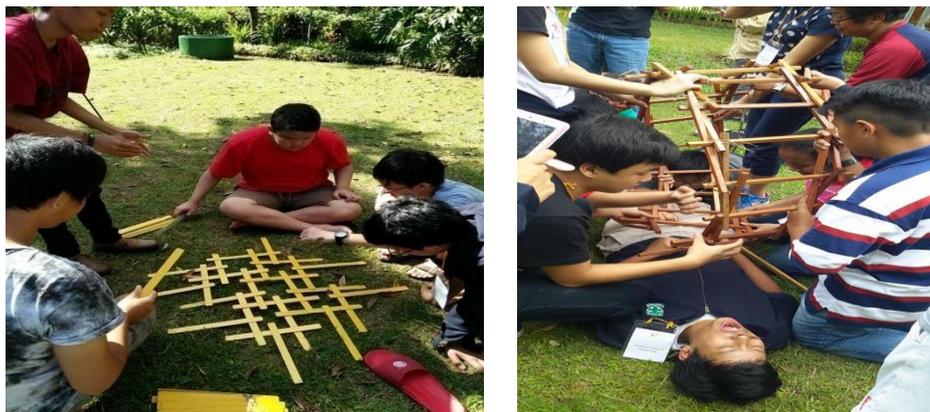

Figure 5. Building Leonardo Dome (left) and Leonardo Bridge (right) (Source :private photo, students are from middle schools of Santa Ursula BSD Jakarta and BPK Penabur Cirebon Indonesia under the project Junior Summer Camp, in 10-19 June 2017,Salatiga, Indonesia June 2017).

### 4.2 Activities in undergraduate students and the related mathematical topics

Fortunately, the same topics for middle schools can also be designed for undergraduate students though the formal mathematical contents must be deeper in this case. Undergraduate students are learning platonic solids and algebraic surfaces in formal manners. Finally, students must find new performances to presents these topics.

### 4.2.1 Parametric curve into curves stitching and surfaces

One of small topics in physics and mathematics is learning parametric curves including the properties. In this research, a modified hypocycloid curve is innovated into many kinds of curves that one may use into ornaments. An old version of the activity is doing curves stitching by representing curves with threads. The same curves can also be ornamented using iron as shown in Fig.6. The variations of curves are obtained due to different values of parameters in the parametric equations. Additionally every curve has its derivative yielding to a new curve. Therefore this activity is also revealing an operator mathematics such as derivative into an object.

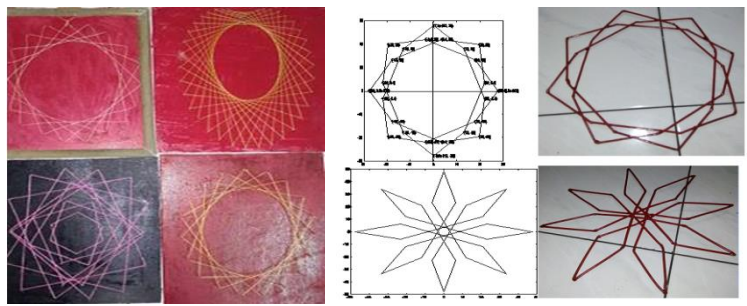

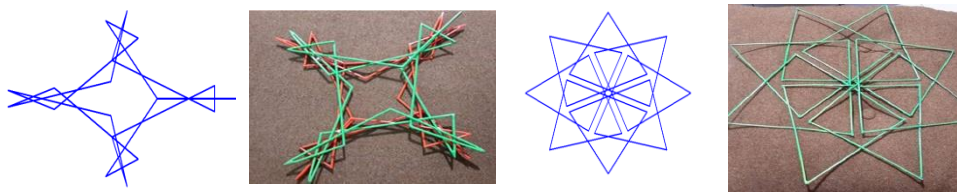

**Figure 6**. A modified hypocycloid curve is defined into several curves stitching and ornaments .

The same parametric equation can also be extended into a surface with spherical coordinates and innovated with other mathematical terms such as Golden ratio, Fibonacci sequence to employ in the equations. Many surfaces can be obtained and some of these are depicted in Fig. 7. Mathematical formulations of surfaces have been presented into detail (Parhusip,2015).

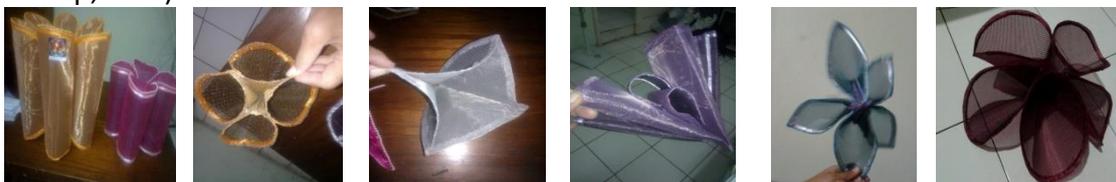

**Figure 7.** Innovative surfaces from a modified hypocycloid curve

The same designs of surfaces are able to be other applications, e.g. motifs for textiles called batik, designs for motifs in bags, souvenirs, and puzzles.

**4.2.2 Geometrical innovations with Rinus Reolofs styles**

Plane geometry is one of topics for schools and undergraduate mathematics students. Innovations have been available in the internet with materials in traditional mathematical questions added with modern presentations such as on line worksheets, using software. These methods sometimes are not known by students and educators in villages in Indonesia. Therefore mobile materials may be suitable choices for students to have interaction between the given topic and an innovation material. The plane geometry is innovated into platonic solids for instance where students may use papers to construct new performance of platonic solids. Other possible material is dealing with outdoor activity by creating Leonardo Dome and Leonardo Bridge as relations between geometrical topic with architecture and physics. Referring to Rinus Reolofs, Leonardo Dome is created using ice sticks and bamboos as the local materials in Indonesia and Leonardo Bridge may be built using woods which are save materials for students and depicted in Fig. 8.

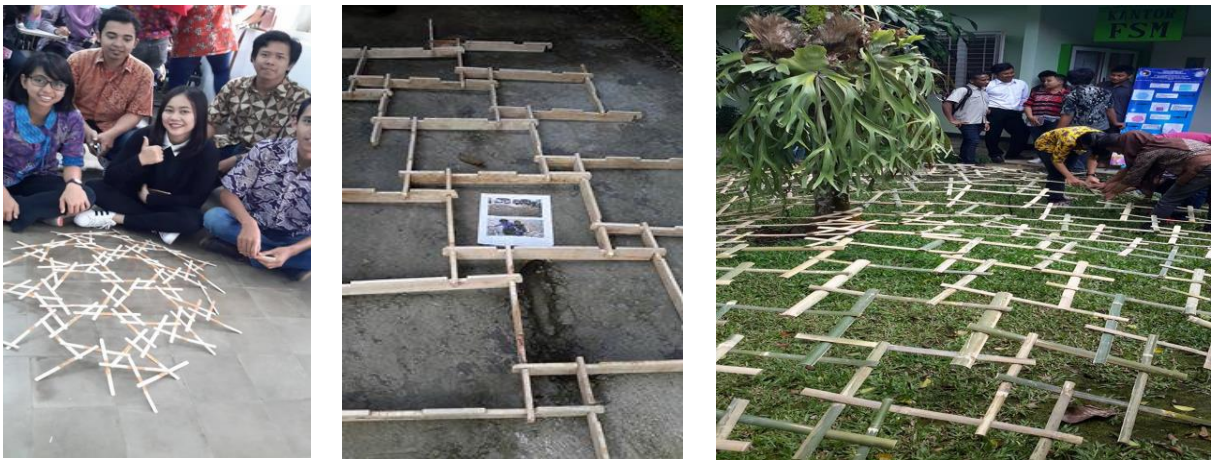

**Figure 8.** Models of Leonardo Dome using stick ices (left), wood (middle) and bamboo (right)

### 4.2.3 Innovation with unused materials

Students may learn triangles and higher degrees of polygon. The pattern can be used for several media using unused material such as bottle aqua. Figure 9 depicts the model of hexagon.

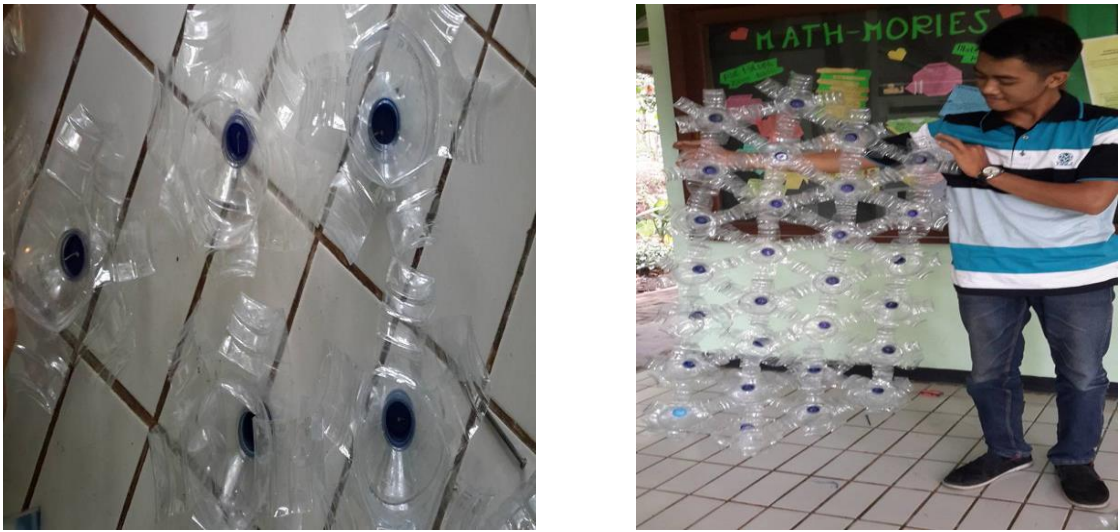

Figure 9. Hexagon patterns are used to create for innovative ODEMA(Ornament Decorative Mathematics)

### 4.2.4 Algebraic Surfaces into Local Cultures in textiles motifs in Batik

Students may also learn the same topic for middle schools such as learning cylinder and sphere as regular algebraic surfaces introduced in classroom. Improving knowledge in algebraic surfaces, students are able to work the surfaces into motifs of textiles called batik as local culture in Indonesia. Using Surfer as friendly software, students create a lot of motifs. Guarded by a mentor from a home industry in Batik, students finally have their own innovation in these algebraic surfaces. Some results are shown in Fig.9.

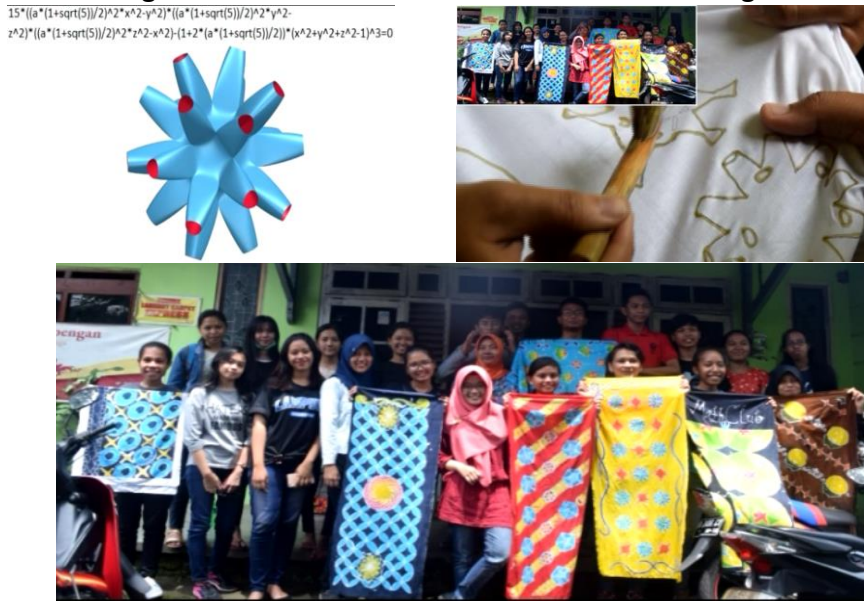

**Figure 9.** Painting textiles (Batik) with algebraic surfaces as motifs from Surfer done by undergraduate mathematics students UKSW, February 2017.

**CONCLUSION**

Several approaches on innovating mathematics are presented here due to research activities done for creating new media for promoting mathematics, geometry particularly. The materials are defined due to the level of students: middle and junior schools and undergraduate students. Platonic solids created from regular planes are innovated into dress or other ornaments following Rinus Reolofs's patterns. In the level of undergraduate students, parametric curves and complex functions are explored into innovative curves and surfaces being new to be new curves for curves stitching, motifs in several kinds of souvenirs and accessories.

Algebraic surfaces are introduced using Surfer for all level students through batik painting activity known to be a usual interesting activity for students offered in schools and universities. This activity becomes a combination between local culture and transferring new technology for students and home industries.